\def\qed{\nobreak\hskip 2pt \vrule height 5pt width 5pt depth 0pt}
 \documentstyle[fullpage,11pt]{article}
 
\begin{document}
\title{ Geometric Finiteness and Uniqueness for  Kleinian
Groups with Circle Packing Limit Sets } \author{Linda
Keen\thanks{
Supported in part by NSF GRANT DMS-8902881} \\Mathematics
Department\\CUNY Lehman College\\Bronx, NY 10468, U.S.A. \and
Bernard Maskit\thanks{Supported in part by NSF GRANT
DMS-9003361}\\Mathematics Department\\SUNY Stony Brook\\Stony Brook, NY
11794-3561, U.S.A. \and Caroline
Series\\Mathematics Institute\\ Warwick University\\Coventry CV4 7AL,
U.K.}
\maketitle
\thispagestyle{empty}
\def\IMSmarkvadjust{0 pt}
\def\IMSmarkhadjust{0 pt}
\def\IMSmarkhpadding{0 pt}
\def\IMSpubltext{Published in modified form:}
\def\SBIMSMark#1#2#3{
 \font\SBF=cmss10 at 10 true pt
 \font\SBI=cmssi10 at 10 true pt
 \setbox0=\hbox{\SBF \hbox to \IMSmarkhpadding{\relax}
                Stony Brook IMS Preprint \##1}
 \setbox2=\hbox to \wd0{\hfil \SBI #2}
 \setbox4=\hbox to \wd0{\hfil \SBI #3}
 \setbox6=\hbox to \wd0{\hss
             \vbox{\hsize=\wd0 \parskip=0pt \baselineskip=10 true pt
                   \copy0 \break%
                   \copy2 \break%
                   \copy4 \break}}
 \dimen0=\ht6   \advance\dimen0 by \vsize \advance\dimen0 by 8 true pt
                \advance\dimen0 by -\pagetotal
	        \advance\dimen0 by \IMSmarkvadjust
 \dimen2=\hsize \advance\dimen2 by .25 true in
	        \advance\dimen2 by \IMSmarkhadjust

%
%
  \openin2=publishd.tex
  \ifeof2\setbox0=\hbox to 0pt{}
  \else 
     \setbox0=\hbox to 3.1 true in{
                \vbox to \ht6{\hsize=3 true in \parskip=0pt  \noindent  
                {\SBI \IMSpubltext}\hfil\break
                {\it J. Reine Angew. Math} {\bf 436} (1993), pp. 209--219.
 
                \vfill}}
  \fi
  \closein2
  \ht0=0pt \dp0=0pt
 \ht6=0pt \dp6=0pt
 \setbox8=\vbox to \dimen0{\vfill \hbox to \dimen2{\copy0 \hss \copy6}}
 \ht8=0pt \dp8=0pt \wd8=0pt
 \copy8
 \message{*** Stony Brook IMS Preprint #1, #2. #3 ***}
}

\SBIMSMark{1991/23}{December 1991}{}
\bibliographystyle{plain}

Let $G \subset PSL(2,{\bf C})$ be a geometrically finite Kleinian
group,
with  region of discontinuity $\Omega (G)$.
By Ahlfors' finiteness theorem, the quotient, $\Omega (G)/G$, is a
finite
union of Riemann surfaces of finite type. Thus on it, there are only
finitely many
mutually disjoint free homotopy classes of simple closed curves.
It is shown in \cite{Maskit7} and \cite{Oh2} that if $\gamma_1,\ldots
,\gamma_k$
are a set of mutually disjoint and simple closed curves
in $\Omega (G)/G$, represented by primitive non-elliptic and
non-conjugate
elements $g_1,\ldots ,g_k$ of $G$, then there is a group $ G'$, and an
isomorphism $\phi \colon G \to G'$ taking parabolic elements of $G$
to parabolic elements of $G'$, for which the images $\phi(g_1),
\ldots, \phi(g_k)$ in $G'$ are parabolic. Heuristically, this means
that the
curves $\gamma_i$ have been ``shrunk'' or ``pinched'' to punctures.
 
In this paper, we assume that $G$ is a finitely generated torsion free
non-elementary Kleinian group with $\Omega(G) \neq \emptyset$. We show
that the maximal number of elements of $G$ that can be
pinched is precisely
 the maximal
number of rank 1 parabolic subgroups that any group isomorphic to
$G$ may contain.
 A group with this largest number of rank 1 maximal parabolic subgroups
 is called {\it maximally parabolic}.  We show such groups exist. We
state our
main theorems         concisely here. Full statements appear
in sections 4 and 5.
 
\smallskip
\noindent
{\bf Theorem I.} {\sl The limit set of a maximally parabolic group
is a circle packing; that is, every component of its regular set is
a round disc.}
 
\noindent and
 
\smallskip
\noindent
{\bf Theorem II.} {\sl A maximally parabolic group is geometrically
finite.}
 
\smallskip
 
We define a class of groups called {\it pinched function groups}.
Roughly speaking, these are either
function groups, or groups isomorphic to function groups,
 for which a set of loxodromic elements has been pinched to
parabolic ones. The precise definition is given in section 5.

 \smallskip
 
\noindent
{\bf Theorem III.} {\sl  A maximally parabolic pinched function group is
 determined up to conjugacy in $PSL(2,{\bf C})$ by its abstract
isomorphism class and its parabolic elements.}

It is a well known consequence of J{\o}rgensen's inequality that
if $G$ is a non-elementary geometrically finite function group, and
 ${\cal T}(G)$ is its quasi-conformal deformation
space, then
any group
in $\partial {\cal T}(G)$ is discrete.
Let $g_1,\ldots ,g_k$ be elements of $G$ which represent a
maximal set of mutually disjoint and non-homotopic simple closed
curves on $\Omega (G)/G$.
An application of our results is that there exists a unique point in
$\partial {\cal T}(G)$ for which all the elements $(g_i), i
= 1,\ldots ,k$, are parabolic.
 
An outstanding conjecture in the theory of Kleinian groups (see for
example, \cite{Th}) is that groups on the boundary of
certain\footnote{A suitable condition
is that the components of the boundary of a compact core of ${\bf
H}^3/G$ are incompressible.}
deformation spaces
of Kleinian groups are
uniquely determined by their ending laminations.  Our result
proves this conjecture in a special case.
 
A well studied deformation space of Kleinian groups is the deformation
space of free groups on two generators. This is the Schottky space of
genus 2; for $G$ in this space $\Omega(G)/G$ is a surface $S$ of genus
2.
The space has three complex dimensions and has a complicated
boundary.
A one dimensional slice of this boundary is formed by shrinking two
curves on $S$, one a dividing curve and the other a non-dividing curve.
The resulting surfaces consist of a punctured torus and a triply
punctured sphere. Since the triply punctured sphere has a unique
conformal structure, this boundary slice is an embedding of the
Teichm\"uller space of a punctured torus. It is known as the Maskit
embedding \cite{Maskit5}. Wright \cite{Wright} studied the boundary of
this embedding
computationally and conjectured that the cusp groups for which an
element corresponding to a particular
simple curve on the torus is parabolic were unique. Our theorem proves
this conjecture.
 
Another one dimensional boundary slice of the Schottky space of genus
2 is formed by shrinking
two non-dividing curves on $S$.  The resulting surfaces are spheres
punctured at four points and the resulting groups are generated by two
parabolics.  Groups generated by two parabolics were investigated
computationally by Riley \cite{Riley} and by Maskit and Swarup
\cite{Mas-Sw};  this boundary slice is known as the
Riley slice. Our theorem proves that the cusp groups on the boundary
of the Riley slice for which
an element corresponding to a particular simple closed curve on the
four punctured sphere is parabolic are unique.
 
On each of these boundaries, three elements of
$G$ were specified to
be parabolic and this determined the group uniquely.  It also follows
from our uniqueness theorem
that if the same three elements are parabolic at a boundary point of a
Maskit slice and at a boundary point of a Riley slice, then  they
correspond
to the same group on the boundary of Schottky space.
  We refer the reader to \cite{KS2,KS4} and the references
therein for further discussion of these spaces.

The method we use for the proofs of theorems I and II is an extension
of that developed by Maskit and Swarup in \cite{Mas-Sw} where they
prove theorem II in the special case that $G$ is generated by two
parabolic elements. The method depends crucially on the fact that any
hyperbolic manifold has a compact core that supports all the homotopy.
  For theorem III we study the deployment of the circles in
the limit set and then apply Marden's isomorphism theorem
\cite{Marden}.
 
The outline of this paper is as follows.
Section 1 contains basic definitions and notation
and a summary of the results we need
about cores
of hyperbolic 3-manifolds. We discuss cusp cylinders and pairing
tubes associated to cyclic parabolic subgroups.
In section 2, we
 define an
 invariant of the isomorphism class of $G$, called the
{\em boundary characteristic} of a Kleinian group,
in terms of the Euler characteristic of the boundary of the core. The
boundary characteristic determines how many loops may be
pinched.
 In section 3 we discuss maximally parabolic
groups and prove that the isomorphism class of any geometrically finite
group contains a group that is maximally parabolic. In section 4 we
prove theorems
I and II. The exact definition of {\em pinched function groups} is
given in section 5 and theorem III is proved there.
 
We would like to thank Jean-Pierre Otal and Colin Rourke for
helpful discussions about this paper. The second author thanks
the Institut des Hautes Etudes Sci\'entifique  and
the last author thanks the M.S.I. at Stony Brook for their
hospitality.

\section{Preliminaries }
\noindent {\bf 1.1. Basics.}
\medskip
 
\noindent
A Kleinian group $G$ is a discrete subgroup of
$PSL(2,{\bf C})$.
It
 acts by conformal automorphisms on the
Riemann sphere ${\hat {\bf C}}$ and by isometries on hyperbolic
three space ${\bf H}^3$.
We denote by $\Omega = \Omega (G) \subset \hat{\bf C}$ the region of
discontinuity
of $G$, {\it i.e.} the set on which the elements of $G$ form a
normal family, and by $\Lambda = \Lambda (G)$ its complement.
Throughout this paper, we make the hypotheses,
which we shall not restate, that all our
Kleinian groups are non-elementary, finitely generated,
torsion free and that $\Omega(G) \neq \emptyset$.
 
A {\it component} of the group $G$ is a connected component of $\Omega$
and a subgroup of $G$ stabilising a component is called a {\it
component subgroup}.
 We denote the quotient ${\bf H}^3 \cup \Omega (G)/G$ by $M = M(G)$. The
group  $G$ is called {\it geometrically finite } if it has a finite sided
convex
fundamental polyhedron in ${\bf H}^3$.

In what follows we shall frequently need to make use of a more
refined notion of invariance called {\it precise invariance}, (see
\cite{Maskit6}). Namely, let $H$ be a subgroup of $G$;
a subset $E \subset \Omega \cup {\bf H}^3$ is said to be {\it
precisely invariant} under $H$ in $G$ if $E$ is invariant under
$H$ and if $g(E) \cap E = \emptyset$ whenever $g \notin H$.
More generally, a family of sets $E_1,\ldots ,E_n$ is precisely
invariant under a family of subgroups $H_1, \ldots ,H_n$ if each
$E_i$ is invariant under
$H_i$ and if, for any $g \in G$, $g(E_i) \cap E_j \ne
\emptyset $ only if $i = j$ and $g \in  H_i=H_j$.

\medskip
 \noindent {\bf 1.2. Parabolic cusps and cores }
\medskip
 
\noindent
 Suppose that $H$ is a maximal parabolic subgroup of rank one in a
 Kleinian group $G$; there is a single point $p$
 fixed by all its elements. Let $C$ be an open horoball ({\it i.e.} an
open
 Euclidean 3-ball) tangent to $\partial {\bf H}^3$ at $p$.
 Clearly $C$ is $H$-invariant and in fact,
 (see {\it e.g.} \cite{Maskit6} p.120), one may
 choose $C$ to be precisely invariant under $H$.
  By putting $p$ at infinity,
 and using the precise invariance, one sees easily that $\partial
 C/G$ is topologically an open annulus if
 $H$ has rank 1, while it is a torus if $H$ has rank 2.

Suppose that $G $ contains $d$
 non-conjugate maximal
parabolic subgroups $ H_1,\ldots , H_d$.
For each $ i = 1,\ldots ,d$, let $ C_i$ be an open
horoball tangent to $\partial
{\bf H}^3$ at the fixed point of $ H_i$.
Such a $ C_i$ is automatically $H_i$ invariant and in fact,
(see \cite{Maskit6} pp. 119-120), one can choose the set $
(C_1,\ldots, C_d)$ to be mutually disjoint and precisely
invariant under $( H_1,\ldots , H_d)$ in $ G $.
It follows from the precise invariance that the $C_i$ all
have disjoint projections, and that the boundaries of
$C_i/G$ are either tori or compact cylinders. We call the
corresponding $C_i/G$ {\it cusp
tori} and {\it cusp cylinders} respectively. Let $\tilde N(G) = {\bf
H}^3 -
 \cup_{i=1}^d C_i$. Then $N(G) =\tilde N(G)/G  =M(G) - \{ \}
{\rm cusp \; tori}\} \cup \{ {\rm cusp \; cylinders} \} \}$.

If $G$ is a Kleinian group containing parabolics, then $M(G)$ certainly
is not compact. A {\it core} $K$ of a 3-manifold $M$ is a compact
connected submanifold of $M$ with the property that the injection $K \to
M$ induces an isomorphism of fundamental groups. It was shown by Scott
\cite{Scott} that every 3-manifold with finitely generated
fundamental group
has a core; this was refined by McCullogh \cite{McCul}, who showed that
the
core can be chosen so that it intersects the boundary of $M$ in any
predetermined compact submanifold of $\partial M$.

If we apply these results to $N(G)$ we see that we can choose a core
$K(G)$ for $N(G)$ so that its boundary contains the cusp tori and cusp
cylinders. Since there exists a retraction from $M(G)$ to $N(G)$, any
core for $N$ is automatically a core for $M$. Notice that because
we are in hyperbolic space the core can be chosen to be irreducible.
 
\medskip
\noindent{\bf 1.3 Pairing tubes}
\medskip
 
\noindent
 When we have a cusp cylinder, in general
we know nothing about the action of $G$ in the
 neighbourhood of a parabolic fixed point outside the
 horoball neighbourhood $C$. However, there is one important case in
 which we have much more precise information. Namely, suppose
 that there are two disjoint circular discs, $D_1,D_2 \subset \Omega $,
 mutually tangent at $p$, whose union is precisely invariant
 under the subgroup $H$ fixing $p$. Then the parabolic point $p$ is
said
 to be {\it doubly cusped}. The two hyperbolic half spaces, $\Sigma_1$
and  $\Sigma_2$
 above $D_1$ and $D_2$, in ${\bf H}^3$ are also precisely
 invariant, and the quotient $\Sigma_1 \cup C \cup \Sigma_2/G$ is a
cylinder
 with two ends at infinity, as described in detail in \cite{Marden}
Sec
 2.6 or \cite{Greenberg} Lemma 4. The quotient cylinder in this
 set-up is called a {\it pairing tube}.
 
\medskip
\noindent
{\bf Definition 1.3.1.}  The {\em natural truncation} $T(G)$ of
$M(G)$ is the manifold obtained by removing the pairing tubes and cusp
tori from $M(G)$.
\medskip
 
 The natural truncation is clearly a retraction of $M(G)$. The
boundary of $T$ is the union of the boundaries of the cusp tori and
components formed from the boundary components
of $\Omega/G$ as follows. We  remove a punctured disc neighbourhood
from each double cusp and
then, for each pair of discs, we glue in an annulus joining the boundary
of
one disc at this cusp to the other. In this way we obtain the boundary
of a cusp cylinder.  Note that the boundary of $T$
is compact if and only if all the rank 1 parabolic points are doubly
cusped.

\medskip
\noindent
{\bf Remark 1.3.2.}  It follows from \cite{Marden} that if $G$ is
geometrically finite then all the rank 1 parabolic fixed points are
doubly cusped. Hence we can restate Marden's criterion for geometric
finiteness in terms of $T(G)$ as: $G$ is geometrically finite if and
only if $T(G)$ is compact.

\section{Boundary characteristic and geometric finiteness }
 
\noindent {\bf Proposition 2.1.} {\sl The maximal number of distinct
disjoint
 simple  closed curves on the boundary of a core of $M(G)$
depends only on the
isomorphism class of $G$.
In particular, it is independent of the choice of core.}

\medskip
\noindent {\bf Proof.}
 Let $G$ and
 $ G^{\prime}$ be isomorphic finitely generated Kleinian groups.
  Let $K(G)$ and
 $K(G^{\prime})$ be cores of $G$ and $G^{\prime}$ chosen as in
section  1.2.
 
 Note that $K(G)$ and $K(G^{\prime})$ are $K(\pi,1)'s$ with isomorphic
 fundamental groups, and hence are homotopy equivalent. In
 particular $K(G)$ and $K(G^{\prime})$ have the same Euler
characteristic.
 Since the Euler characteristic of a closed orientable
 3-manifold is zero, the Euler characteristic of the boundary of
 a compact orientable 3-manifold is twice the Euler
 characteristic of the manifold itself.
 Hence the boundaries of $K(G)$ and $K(G^{\prime})$ have the same
Euler
 characteristic, say $2 \chi$.

As remarked above, we may assume  neither core has a boundary
component that is
 a sphere.
 Since the groups $G$ and $G^{\prime}$ are isomorphic, and the tori on
the
 boundary of the core correspond exactly to the conjugacy
 classes of maximal rank 2
 parabolic subgroups (equivalently rank 2 abelian),
$K(G)$ and $K(G^{\prime})$
 have the same number
 $\tau$ of tori on their boundaries.
 
The maximal number of non-homotopic disjoint simple closed
 curves on a union of surfaces of negative Euler characteristic
 is exactly $-3 \chi$. Hence, for both
$K(G)$ and $K(G^{\prime})$,
 the
 number of disjoint homotopically distinct simple loops on their
 boundaries is $\tau - 3\chi$.\qed
 
 \smallskip
 \noindent{\bf Definition 2.2.}
 The {\em boundary characteristic} $b = b(G)$ of the
 Kleinian group $G$  is $ n-\tau$ where $n$ is the maximal number of
 disjoint homotopically distinct simple closed curves on the
 boundary of any core for $G$, and $\tau$ is the number of
 conjugacy classes of rank 2 maximal parabolic subgroups of $G$.
 
\smallskip
We can give a criterion for geometric finiteness of $G$ in
terms of the boundary characteristic and the truncation $T(G)$ of
$M(G)$.
 
\smallskip
\noindent{\bf Proposition 2.3.} {\sl Let $G$ be a Kleinian group. If $G$ is
geometrically finite then the non-toral components of $\partial T(G)$
support $b(G)$
 homotopically distinct disjoint simple
closed curves. The converse is also true if, in addition, all
rank one parabolic fixed points are doubly cusped.}
 
\medskip
\noindent{\bf Proof.} If $G$ is geometrically finite, then by remark
1.3.2
 the natural truncation $T$ is
compact.
Since $T$ is a retraction of $M$, we may take it to be a core
$K$ of $M$; hence it supports at most $b$ distinct simple closed
curves on the non-toral components.
 
Suppose conversely that the non-toral part of $\partial T$ supports $b$
distinct simple closed curves, and that every rank one parabolic fixed point
is doubly cusped. Let $\Omega_1$ be the set obtained from $\Omega$ by
removing disc neighbourhoods of the doubly cusped parabolic points.  The
boundary $\partial T$ is the union of the boundaries of the pairing
tubes, the cusp tori and the components of $\Omega_1/G$. Using  the
assumption about the rank one parabolic points, we see that $\partial T$ is
compact.  Thus we may choose
a compact core $K$ for $M$ whose boundary contains $\partial T$ and
such that $K \subset T$.
 By proposition 2.1, the boundary
$\partial K$ of
$K$ also supports $b$ disjoint simple closed curves on its non-toral
part. Now $\partial K$ and $\partial T$ contain the same number of cusp
tori so $\partial T = \partial K$ and hence $K=T$.
Hence $T$ is compact and so Marden's
criterion for
 geometric finiteness is satisfied (see \cite{Marden} and also
\cite{Maskit6} pp. 128-130).\qed
 
\smallskip
\noindent{\bf Remark.}
Observe that if $G$ is geometrically finite and purely
 loxodromic,  then $T= M$, $\partial T = \partial M$ and $b$
is the maximal
 number of disjoint homotopically distinct simple closed curves
 that one can draw on $\Omega/G$.
 
\section{Boundary characteristic and parabolics}
 
\noindent{\bf Proposition 3.1.} {\sl The number of conjugacy classes
of rank 1 maximal parabolic subgroups in a Kleinian group $G$ is
at most $b(G)$.}
 
 \smallskip
 \noindent{\bf Proof.} Form the core $K$ as described in section
 1.2. For each cusp cylinder in $\partial{K}$ we have a simple loop,
and these loops are disjoint and homotopically distinct. Since there are at
most
$b$ such cylinders we are done. \qed
 
\smallskip
\noindent{\bf Definition 3.2}
 A Kleinian group $G$ that contains
 $b(G)$ distinct conjugacy classes of rank 1 maximal parabolic
 subgroups is called  {\em maximally parabolic}.
 
As an immediate corollary to proposition 3.1 we have:
 
 \noindent{\bf Corollary 3.3.} {\sl
  If $G$ is maximally parabolic, then, on any core $K(G)$
there are $b(G)$ disjoint homotopically distinct simple
 closed curves on non-toral components of $\partial K(G)$, each of
which
  is represented by a parabolic element of
 $G$. These $b(G)$ simple loops divide the non-toral part of $\partial
K(G)$ into
 $2b(G)/3$ pairs of pants.}
 
 \smallskip
The existence of maximally parabolic groups follows from the following
theorem.
\smallskip
 
 \noindent {\bf Theorem 3.4.} {\sl Let $G$ be a
geometrically finite Kleinian group. Then there exists a maximally
parabolic group $G'$ and an isomorphism $\phi : G \to  G'$
 mapping parabolic elements of $G$ to parabolic
 elements of $ G'$.
 }
 
\smallskip
The essential ingredient in the proof is the following theorem.
\smallskip
 
\noindent{\bf Theorem 3.5.} {\sl Let $G$ be a geometrically
finite
Kleinian group and let $g \in G$ be a loxodromic element representing
a simple closed loop on some component of $\Omega/G$. Then there
exists a geometrically finite Kleinian group $G'$, and an isomorphism
$\phi \colon G \to G'$ that maps parabolic elements of $G$ to
parabolic elements of $G'$, and such that   $\phi(g)$ is  parabolic.}

\smallskip
This theorem was proved by Maskit
\cite{Maskit7} for function groups and the proof is  easily extended to
the class of pinched function groups defined in section 5.
Based on an extension of
theorem 7.1 in \cite{Th2},
 Ohshika, \cite{Oh2}, generalised this result
 to all geometrically
finite Kleinian groups.

\smallskip
We shall also need the following lemma.
 
\smallskip
\noindent{\bf Lemma 3.6.} {\sl Let $G$ be a Kleinian group and
suppose
that some component $S$ of $\Omega/G$ is not a triply punctured sphere.
Then there is a loxodromic element in $G$ that represents a simple
closed curve on $S$.}
 
\medskip
\noindent{\bf Remark.} This lemma is closely related to proposition 2
in \cite{Oh2}.
 
\medskip
\noindent{\bf Proof.} Let $\Omega_0$ be a component of $\Omega$ that is
a lift of the surface $S$.  The component subgroup
$H$ stabilising $\Omega_0$ is non-elementary since its limit set
$\partial\Omega_0$ is a subset of the perfect set $\Lambda(G)$. Thus,
$H$ is a function group with invariant component $\Omega_0$ containing
infinitely many loxodromic elements, all of which represent
loops on $S$. Assume that none of these is
simple.
 
Using combination theorem techniques, it was shown in \cite{Maskit8}
(see also \cite{Maskit6}, chapter X) that
$H$ is a free product of a finite
number of subgroups $H_1,\ldots, H_n$, where each factor is one of
three possible types. The factor may be a group with a
simply connected invariant component;
in this case, we have a universal cover for a subsurface of $S$. It
may be an elementary group with one limit point; then it is
either a rank 1 or rank 2 parabolic subgroup. Finally, it may be a
purely loxodromic free group ({\em i.e.} a Schottky group).
 
One of the conclusions of the free product combination theorem is that
one may obtain a fundamental domain for $H$ from fundamental domains for
the factors $H_i$ by cutting and pasting fundamental domains for the
$H_i$ in such a way that if $h \in H_i$ represents a simple loop on
$\Omega(H_i)/H_i$, then $h$ represents a simple loop on $S$. We may
therefore assume that no loxodromic element of any $H_i$ represents a
simple loop.
 
The quotient surface of a purely loxodromic free group
 is
a compact surface and each generator represents a simple loop. Hence, by
our assumption, none of the factors $H_i$ are of this type.

It was shown in
\cite{Maskit8} that
if $H_i$ has a simply connected invariant domain, then the set of
accidental parabolic elements represents a finite set of simple disjoint
homotopically distinct loops on its quotient surface.
Either there
are infinitely many simple loops and so most are
 represented by  loxodromic elements or
 $H_i$ is a Fuchsian
group with
 quotient a triply punctured sphere.
 
Since $H$ is non-elementary, and itself not a triply punctured sphere
group, we conclude that the free product description of
$H$
has at least two factors and that
 the factors are either distinct rank 1 or rank 2 parabolic
subgroups or triply punctured sphere groups. In the cut and paste
construction
of fundamental domains, for each pair of factor subgroups, there is a
homotopically nontrivial loop, say $\tilde \gamma$,
so that fundamental domains for this pair of
factors are pasted across $\tilde{\gamma}$. This loop projects to a
simple
homotopically non-trivial loop $\gamma$ on $S$. Since $S$ is hyperbolic
of finite
area, and not a triply punctured sphere, there is a simple homotopically
non-trivial loop on $S$ crossing $\gamma$ that is represented by a word
of length at least two (using the free product word length) in $H$.
Another conclusion of the free product combination theorem is that every
parabolic element in $G$ is conjugate to a parabolic element of some
factor.  It
follows, therefore, that the element we have found is
loxodromic. \qed

\smallskip
\noindent{\bf Proof of Theorem 3.4}  Let $G$ be a geometrically
finite Kleinian group that is not maximally parabolic. Suppose that it
contains $p(G) < b(G)$ distinct conjugacy classes of rank 1
parabolic
subgroups.  We shall show that there is a group $G'$ and an isomorphism
$\phi \colon G \to G'$, taking parabolics to parabolics, such that $G'$
has $p(G') > p(G)$ distinct conjugacy classes of rank 1 parabolic
subgroups. Since we can repeat this argument until the number of such
classes is equal to $b(G)$, this will prove the proposition.

Since $G$ is geometrically finite, by remark 1.3.2 the
truncation $T(G)$ is compact and so can be taken as a core $K(G)$
for $M(G)$. By proposition 2.3 we can find exactly $b(G)$ disjoint
simple closed loops on the non-toral part of $\partial T(G)$; $p(G)$
of these can be taken to be loops
around the boundaries of the cusp cylinders. If all the components of
$\Omega/G$ were triply punctured spheres it would be impossible to
find any more disjoint loops on $\partial T(G)$ so that
$p(G)=b(G)$; this is not the case.
 
Let $S$ be a component of $\Omega/G$ that is not a triply punctured
sphere. By lemma 3.6 we can find a loxodromic element $g \in G$
representing a simple closed curve on $S$ and by theorem 3.5 there
is a geometrically finite group $G'$ and an isomorphism $\phi \colon G
\to G'$, taking parabolics to parabolics and such that $\phi(g)$ is also
parabolic. Therefore $G'$ has at least one more conjugacy class of rank
1 parabolic subgroups than $G$ as required.
 \qed
 
\section{Proofs of theorems I and II}
 
 In this section we prove the following two related
results. We recall that $G$ is non-elementary, finitely generated,
torsion free and has a non-empty regular set.

\smallskip
 \noindent{\bf Theorem I.} {\sl Let $G$ be a
maximally parabolic Kleinian group $G$ of boundary characteristic
$b(G)$. Then $\Omega(G)$ is a union of round discs and $\Omega(G)/G$ is
a union of triply punctured spheres.
For
each maximal parabolic subgroup $H \subset
 G$ of rank 1, exactly two of these round discs are
 tangent at the fixed point of $H$; in particular, all parabolic fixed
points are doubly cusped.}
 
\medskip
 
\noindent{\bf Theorem II.} {\sl A maximally parabolic Kleinian group
 $G$ is geometrically finite.}
 
An immediate corollary is:
 
 \smallskip
 \noindent {\bf Corollary 4.1.} {\sl Let $G$ be a geometrically
finite Kleinian group
 of boundary characteristic $b(G)$. Let $\phi$ be an
 isomorphism onto a subgroup $\phi(G)$ of $PSL(2,{\bf C})$, where $\phi
 (h)$ is parabolic for every parabolic $h \in  G$, and $\phi(G)$ has
 $b(G)$ distinct conjugacy classes of rank 1 maximal parabolic
 subgroups. Then $\phi(G)$ is either not discrete, or it is discrete
 with non-empty regular set and is geometrically finite.}
 
\smallskip
 \noindent {\bf Proof of theorem I.}
Let $b$ be the boundary characteristic of $G$,
let $K$ be a core for $M$ as in section 1.2, and let $C$ denote the
union of components of $\partial K$ that are not cusp tori. By
corollary 3.3 there are b disjoint homotopically distinct simple
loops on $C$ representing distinct parabolic elements of $G$. These
divide $C$ into 2b/3 pairs of pants.
 
Let $P$ be one of these pairs of pants and let
$\gamma_1,\gamma_2,\gamma_3$ be loops going around the boundary curves
of $P$, chosen so that $\gamma_1 \gamma_2 \gamma_3 = id$ in $\pi_1(P)$.
Each $\gamma_i$ lies on the boundary of a cusp cylinder corresponding to
a conjugacy class of maximal parabolic subgroups in $G$; hence each
$\gamma_i$ bounds a punctured disc contained in the cusp cylinder.
Adjoining these punctured discs to $P$ we obtain a triply punctured
sphere $Q$.
 
Lift $Q$ to ${\bf H}^3$ and let $\tilde Q$ be a connected component of
this lift. Let $R \subset \tilde Q$ be a simply connected fundamental
domain for the action of $G$ on $\tilde Q$, such that the Euclidean
closure of $R$ in $ \hat{\bf C}$ consists of four parabolic
fixed points $\xi_1,\xi_2,\xi_3,\xi_4$ on $\hat{\bf C}$. We label these so
that there is a connected component of the lift $\tilde
\gamma_i$ of $\gamma_i$ with both endpoints at $\xi_i$, $i=1,2,3$. We
can choose generators, $h_i$ of the parabolic subgroups
stabilising $\xi_i$, such that
 $h_1 h_2 h_3 = 1$ in $G$.
An easy computation with matrices now
shows that $F=<h_1,h_2,h_3>$ is Fuchsian, and that its limit set is a circle
$\Lambda(F)$.
 
By comparing the action of $F$ on $R$ with its action on a standard
fundamental domain in the Poincar\'e disc, we see that $\cup_{f \in F}
f(R)$ is both open and closed in $\tilde Q$. Therefore,
$\tilde
Q = \cup_{f \in F} f(R)$. It follows that the boundary of $\tilde Q$ in
${\bf H}^3 \cup \hat{\bf C}$ is exactly $\Lambda(F)$.
 
Since $Q=R/G = \tilde Q/F$ is a triply punctured sphere, we see
that $F
= \pi_1(\tilde Q/F)$ . Consequently, $\tilde Q$ is a universal
cover
for $Q$. It follows that $\tilde Q$ is a topological disc in ${\bf H}^3$
whose boundary in ${\bf H}^3 \cup \hat{\bf C}$ is exactly $\Lambda(F)$.
Thus $\tilde Q$ separates ${\bf H}^3$.
 
Now let $\tilde K$ be the lift of the core $K$ to ${\bf H}^3$. Since $K$
is connected so is $\tilde K$; also, since $Q \subset \partial K$, we know
that $\tilde Q \subset \partial \tilde K$. It follows that $\tilde K$ is
entirely contained on one side of $\tilde Q$.  Denote by $\Delta(F)$ the
open disc in $\hat{\bf C}$ bounded by $\Lambda(F)$ and on the opposite
side of $\tilde Q$ from $\tilde K$. Because $\tilde K$ is $G$-invariant,
there can be no limit points of $G$ in $\Delta(F)$.
 
We next show that $\Delta(F)$ is precisely invariant under $F$. Clearly,
any image of $\Delta(F)$ that intersects $\Delta(F)$ coincides with it,
so we have only to show that the component subgroup
$F_{\Delta}$ of $\Delta$ is just $F$.
Temporarily think of $\Delta(F)$ as a hyperbolic disc in which
$F_{\Delta}$ and $F$ act as Fuchsian groups. If $F_{\Delta} \supset F$,
$F_{\Delta} \neq F$, then $\Delta(F)/F \to \Delta(F)/F_{\Delta}$
would
be a covering of hyperbolic surfaces of degree greater than $1$, and so
$\Delta(F)/F_{\Delta}$ would be a surface
of hyperbolic area less than $2 \pi$. Since $G$, and hence $F_{\Delta}$,
contain no elliptics, this is impossible.

Suppose now that $P$ and $P^{\prime}$ are distinct pairs of pants
adjacent along the same curve on $\partial K$.
Then the circles $\Lambda(F)$ and $\Lambda(F^{\prime})$ constructed
above are
distinct and are both invariant under a common parabolic element $h$  of
$G$,
so they are mutually tangent.  Since
the disks $\Delta$ and $\Delta^{\prime}$ do not
intersect,
 the parabolic subgroup generated by $h$ is doubly cusped in
the sense of section 1.3.  It is clear that the fixed point of every
cyclic parabolic subgroup is such a point of mutual tangency and so
doubly cusped.
 
Finally, we want to show that the union of the discs $\Delta_i$ for each
of the component subgroups $F_i$ is all of $\Omega$; i.e. $\Omega =
\cup_i \Delta_i$. Adjoin to $K$ the portions of the
cusp cylinder bounded between the punctured discs that we adjoined to
the pants in $\partial K$ to make the triply punctured spheres $Q_i$,
and call the resulting region $A$. Then the lift $\tilde A$ of $A$ to
${\bf H}^3$ is connected and its boundary in ${\bf H}^3$ is exactly the
union of all the lifts $\tilde Q_i$ of the triply punctured spheres
$Q_i$.

Choose $\eta \in \Omega $, and pick an oriented geodesic $\beta
$ running from some point in $\tilde A$ to $\eta$.
Let $\xi $ be the last point of $\beta $ in $\tilde A$.
Then $\xi $ lies in
 one of the surfaces $Q_i$.
Since by the choice of $\xi$  the segment of $\beta $ from $\xi $ to
$\eta $ does not recross $\partial \tilde A$, it is easy to see
that $\eta $ lies in a disc $\Delta_i$ as required.\qed
 
\smallskip
 \noindent {\bf Proof of Theorem II.}
By Theorem I we know that all the rank 1 parabolic fixed points are
doubly cusped. It follows that there are the maximal number, $b$, of
simple closed loops on $\partial T$, and hence, by the
finiteness criterion of proposition 2.3, that $G$ is geometrically
finite. \qed

\section{Pinched function groups and uniqueness}
 
We shall prove our uniqueness result for a restricted class of groups
that we call {\it pinched function groups}.
 
Let $G$ be a Kleinian group and let $\{\Omega_i\}$ be the components of
$\Omega$. We say $\Omega_i$ is {\em adjacent} to $\Omega_j$
if there is a doubly cusped parabolic element in $G$ with
invariant discs $D \subset \Omega_i$ and $D' \subset \Omega_j$. We
define an equivalence relation on $\{\Omega_i\}$ by $\Omega_i \equiv
\Omega_j$ if and only if there is a finite sequence of
components $\Omega_{i}=\Omega_{n_1}, \Omega_{n_2}, \ldots
\Omega_{n_d}=\Omega_j$ with $\Omega_{n_k}$ adjacent to
$\Omega_{n_{k+1}}$, $k=1, \ldots, d-1$.
 
We call the union of all the $\Omega_i$ in an equivalence class an {\it
augmented component} of $G$.
 
\smallskip
\noindent{\bf Definition 5.1.} A Kleinian group $G$ is a {\em pinched
function group} if it has a $G$-invariant augmented component.

\smallskip
\noindent{\bf Remark.} The augmented components of a Kleinian group
$G$ are in  one to one
correspondence with the non-toral components of $\partial T$, the
boundary of the natural truncation $T$ of $M$.
 
 \smallskip
It follows that if $G$ is a geometrically finite pinched function group and
if $G^{\prime}$ is
obtained from $G$ by pinching, then $G^{\prime}$ is also a pinched function
group.

\medskip
 
We can now state our uniqueness theorem.
 
\smallskip
 \noindent{\bf Theorem III.} {\sl Let $G$ and $G^{\prime}$ be
 maximally parabolic  Kleinian groups where $G$ is a pinched
function group.
 Suppose that there is an isomorphism $\phi : G \to G^{\prime}$
where
 $\phi(g)$ is parabolic if and only if $g$ is parabolic. Then
 there is a (possibly orientation reversing) M\"obius
 transformation $\alpha$, so that $\phi(g) = \alpha\circ
 g\circ\alpha^{-1}$, for all $g \in G$.}
 
\smallskip
Before we can proceed with the proof we need some preliminaries.  We suppose
that $G$, $G^{\prime}$ and $\phi$ are as in the statement of Theorem III.
 
  Let ${\cal F}(G)$ denote the set of Fuchsian
 subgroups of $G$ with three parabolic generators whose product
 is the identity. Any group $H$ in $\cal F(G)$ has a circular limit
 set contained in $\Lambda (G)$.
 However, as can be seen from
 Example VIII.G.1 in \cite{Maskit6} there are, at least for some
groups
 $G$, elements of ${\cal F}(G)$ whose limit circles have limit
 points of $G$ both inside and outside. We denote the subset of ${\cal
 F}(G)$ consisting of those subgroups which are also component
 subgroups of $G$ by ${\cal F}_0(G)$. We call such groups {\it
peripheral}.
 
 Note that if $H \in {\cal F}_0(G)$, then the the limit set of $H$ is a
circle bounding an
 open disc $\Delta(H)$ in $\Omega(G)$.
 
 It is clear from the definition that an isomorphism
$\phi$ such as the one in theorem III maps
 the groups in ${\cal F}(G)$ onto the groups in ${\cal F}
(G')$.
 
\smallskip
 \noindent{\bf Proposition 5.2.}
{\sl If $G$ is a pinched function group then the isomorphism $\phi$ maps
the groups in ${\cal
 F}_0(G)$ onto those in ${\cal F}_0(G^{\prime})$.}
 
The main ingredient in the proof of this proposition is the
following theorem of Susskind \cite{Suss}:

 \noindent{\bf Theorem 5.3 (Susskind).} {\sl Let $H_1$ and $H_2$ be
 geometrically finite subgroups of a (torsion free) Kleinian group $G$.
  Then
$$\Lambda(H_1) \cap \Lambda(H_2)= \Lambda(H_1\cap H_2).$$}
 
\smallskip
An immediate corollary of theorem 5.3 is:
 
\smallskip
\noindent{\bf Corollary 5.4.} {\sl Let $G$ be a maximally parabolic
Kleinian group and let $F_1,F_2 \in {\cal F}_0(G)$. Then the
circles $\Lambda(\phi(F_1))$ and $\Lambda(\phi(F_2))$ are either tangent
or disjoint. In the first case, $F_1 \cap F_2$ consists of a rank 1
parabolic subgroup; in the second $F_1 \cap F_2 = \{id\}$.}
 
\smallskip
\noindent{\bf Proof.} Since the groups $F_1$ and $F_2$ are peripheral,
by theorem I and Susskind's theorem either $\Lambda(F_1) \cap
\Lambda(F_2) = \emptyset$
and $F_1 \cap F_2 = \{ id \}$, or $F_1 \cap F_2$ is a
rank one parabolic subgroup and
$\Lambda(F_1)$ and $\Lambda(F_2)$ are tangent. Thus, either $\phi(F_1)
\cap \phi(F_2) = \{ id \}$, or
$\phi(F_1) \cap
\phi(F_2)$ is a rank 1 parabolic subgroup. Using Susskind's theorem
again,
$\Lambda(\phi(F_1)) \cap \Lambda(\phi(F_2))$ is either empty or a single
point. \qed

\smallskip
 
 \noindent{\bf Proof of Proposition 5.2.} Suppose that $F_0 \in
{\cal F}_0(G)$ and $\phi(F_0) \not\in {\cal F}_0(G')$; that is,
$\phi(F_0)$
is not peripheral. Then there are points of $\Lambda(G')$ both inside
and outside $\Lambda(\phi(F_0))$.

Let $\Sigma(G)$ denote the invariant augmented component of $G$. Clearly
$\partial \Sigma(G) = \Lambda(G)$.  Now $\Sigma(G)$ is a union of discs
$\Delta_i$ for some collection of $F_i \in {\cal F}_0(G)$. Set $J= \{F_i
\colon \Delta(F_i) \subset \Sigma(G) \}$ and let $B_{\phi} =
\cup_{F_i \in J} \Lambda(\phi(F_i))$.
 
By corollary 5.4,
$B_{\phi}$ is a union of mutually tangent circles.  Since $\Sigma(G)$ is
invariant, it follows  that
 $\bar B_{\phi} = \Lambda(G')$.

Thus we can find $F,F' \in J$ and $\xi \in \Lambda(\phi(F)),
\xi' \in \Lambda(\phi(F'))$ with $\xi$ inside $\Lambda(\phi(F_0))$
and
$\xi'$ outside $\Lambda(\phi(F_0))$.  Now since both $F$ and $F'$ are
in $J$, we have $\Delta(F) \equiv \Delta(F')$ and there is a finite
chain of groups $F=F_1, F_2, \ldots , F_k=F'$  $F_i \in J$, such that
$\Lambda(F_i)$
and $\Lambda(F_{i+1})$ are mutually tangent, $i=1,\ldots,k-1$.
 
Again applying corollary 5.4, the circles $\Lambda(\phi(F_i))$ and
$\Lambda(\phi(F_{i+1}))$ are also mutually tangent, $i=1,\ldots,k-1$,
and furthermore none of the $\Lambda(\phi(F_i))$ can intersect
$\Lambda(\phi(F_0))$ in more than one point. Hence we can find a triple
of circles $\Lambda(\phi(F_0)), \Lambda(\phi(F_r)),
\Lambda(\phi(F_{r+1}))$ with a common point of tangency. Susskind's
theorem then tells us that the circles $\Lambda(F_0),\Lambda(F_r)$ and
$\Lambda(F_{r+1})$ must also have a common point of tangency, but this
is impossible since the three groups $F_0, F_r$ and $F_{r+1}$ are all in
${\cal F}_0(G)$. \qed
 
\medskip
\noindent{\bf Proof of Theorem III.}
Let $F \in {\cal F}_0(G)$ and let $\Delta(F)$ be its invariant disc in
$\Omega(G)$.  By proposition 5.2, $\phi(F) \in {\cal F}_0(G')$ and so
we can find an associated disc $\Delta(\phi(F)) \subset \Omega(G')$.
There is clearly a unique conformal (or anti-conformal) homeomorphism
$\Phi$ of $\Delta(F)$
onto $\Delta(\phi(F))$ compatible with the actions of $F$ and $\phi(F)$.
Since $\Delta(F)$ and $\Delta(\phi(F))$ are precisely invariant under
$F$ and $\phi(F)$ in $G$ and $G'$ respectively, $\Phi$
extends equivariantly to all images of $\Delta(F)$ under $G$. Repeating
the argument for each $F \in {\cal F}_0(G)$ and using the fact, proved
in theorem I, that $\Omega(G) = \cup_{F \in {\cal F}_0(G)} \Delta(F)$,
we obtain a map $\Phi$ of $\Omega(G)$ into $\Omega(G')$
inducing $\phi$.

We now want to apply
Marden's
isomorphism theorem (\cite{Marden} theorem 8.1) to conclude that $\phi$
is induced by a conjugacy
in $PSL(2,{\bf C})$.  By theorem II, $G$ is geometrically finite. We
have not checked the
hypothesis in the isomorphism theorem which requires
that the induced homeomorphism have the same orientation on each
component; however, as observed by Marden and Maskit, \cite{Mar-Mas},
the proof does not use this hypothesis in an
essential way. (See also \cite{Hemp}, theorem 13.9.)  All that is
necessary
is that $\Phi$ have the same orientation type on adjacent components of
$\Omega$. If $\Delta(F)$ and $\Delta(F')$ are adjacent then $F$ and
$F'$ contain a common parabolic which forces the orientation type to be
the same in both components.
 We may therefore apply Marden's theorem to
conclude that $\Phi$ is a conformal (or anticonformal) homeomorphism of
$\hat{\bf C}$.
 \qed

\bibliography{ref}
 \end{document}